\documentclass[11pt]{article}
\newenvironment{proof}{\noindent {\em Proof}.\ }{\hspace*{\fill}$\halmos$\medskip}

\newcommand{\halmos}{\rule{1ex}{1.4ex}}
\newcommand{\QED}{\hfill \halmos} %put !#at right margin
\usepackage{amsmath, amssymb}
\usepackage{graphicx}

\newtheorem{thm}{Theorem}%[section]
\newtheorem{lemma}{Lemma}%[section]
\newtheorem{define}{Definition}%[section]
\newtheorem{cor}{Corollary}%[section]
\newtheorem{pro}{Proposition}%[section]
\newtheorem{rk}{Remark}

\newcommand{\arrowschem}[2]{\raisebox{-2ex}%
	{$\stackrel{\stackrel{\displaystyle#1}{\longrightarrow}}%
	{\stackrel{\longleftarrow}{#2}}$}}

\begin{document}

\title{Singularly Perturbed Monotone Systems and an Application to Double Phosphorylation Cycles}

\author{Liming~Wang
       and~Eduardo~Sontag%
\thanks{L. Wang is with the Rutgers University, Department of Mathematics,
\texttt{wshwlm@math.rutgers.edu}.}
\thanks{E. Sontag is with the Rutgers University, Department of Mathematics,
\texttt{sontag@math.rutgers.edu}.}
}

\maketitle

\begin{abstract}
The theory of monotone dynamical systems has been found very useful in the
modeling of some gene, protein, and signaling networks.
In monotone systems, every net feedback loop is positive.
On the other hand, negative feedback loops are important features of many
systems, since they are required for adaptation and precision.
This paper shows that, provided that these negative loops act at a
comparatively fast time scale, the main dynamical property of (strongly)
monotone systems, convergence to steady states, is still valid.
An application is worked out to a double-phosphorylation ``futile cycle''
motif which plays a central role in eukaryotic cell signaling.
\end{abstract}

\section{Introduction}
Monotone dynamical systems constitute a rich class of models, for which global
and almost-global convergence properties can be established.  They are
particularly useful in biochemical applications and also appear in areas like
coordination \cite{moreau} and other problems in control \cite{chisci}.
One of the fundamental results in monotone systems theory is Hirsch's Generic
Convergence Theorem \cite{Hirsch2,Hirsch,Hirsch-Smith,Smith}.
Informally stated, Hirsch's result says that almost every bounded solution of
a strongly monotone system converges to the set of equilibria.  There is a rich
literature regarding the application of this powerful theorem, as well as of
other results dealing with everywhere convergence when equilibria are unique
(\cite{Dancer,JiFa,Smith}), to models of biochemical systems.
See for instance~\cite{04sysbio,05ejc} for expositions and many references.

Unfortunately, many models in biology are not monotone, at least with respect
to any standard orthant order.  This is because in monotone systems (with
respect to orthant orders) every net feedback loop should be positive, but, on
the other hand, in many systems negative feedback loops often appear as well,
as they are required for adaptation and precision.
However, intuitively, negative loops that act at a
comparatively fast time scale should not affect the main characteristics
of monotone behavior.
The main purpose of this paper is to show that this is indeed the case, in the
sense that singularly perturbed strongly monotone systems
inherit generic convergence properties.
A system that is not monotone may become monotone once that fast variables are
replaced by their steady-state values.
In order to prove a precise time-separation result, we employ tools
from geometric singular perturbation theory.

\begin{figure}[h]
  \centering \includegraphics[scale=0.4,angle=0]{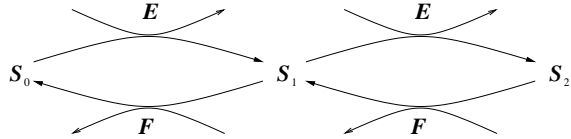}
  \caption{Dual futile cycle. 
A substrate $S_0$ is ultimately converted
into a product $P$, in an ``activation'' reaction triggered or facilitated by
an enzyme $E$, and, conversely, $P$ is transformed back (or ``deactivated'')
into the original $S_0$, helped on by the action of a second enzyme $F$.}
\label{fig:Scheme}
\end{figure}

This point of view is of special interest in the 
context of biochemical systems; for example, Michaelis Menten kinetics are
mathematically justified as singularly perturbed versions of mass action
kinetics~\cite{keshet,murray-book}.
One particular example of great interest in view of current systems biology
research is that of ``futile cycle'' motifs, as
illustrated in Figure \ref{fig:Scheme}.
As discussed in~\cite{Samoilov}, futile cycles (with any number of
intermediate steps, and also called substrate cycles, enzymatic
cycles, or enzymatic interconversions) underlie signaling processes
such as GTPase cycles \cite{Donovan}, 
bacterial two-component systems and phosphorelays 
\cite{groisman, grossman}
actin treadmilling \cite{chen}),
and glucose mobilization \cite{karp},
as well as metabolic control \cite{stryer}
and cell division and apoptosis \cite{sulis} and cell-cycle checkpoint control
\cite{lew}.
A most important instance is that of Mitogen-Activated Protein Kinase
(``MAPK'') cascades, which regulate primary cellular activities such as
proliferation, differentiation, and apoptosis 
\cite{lauffenburger,Chang,
ferrell,widman}
in eukaryotes from yeast to humans.
MAPK cascades usually consist of three tiers of similar structures with
multiple feedbacks \cite{Burack}, \cite{Ferrell_Bhatt}, \cite{Zhao}. Here we
focus on one individual level of a MAPK cascade, which is a futile cycle
as depicted in Figure
\ref{fig:Scheme}.
The precise mathematical model is described later.
Numerical simulations of this model have suggested that the system may be
monostable or bistable, see \cite{Kholodenko}. The later will give rise to
switch-like behavior, which is ubiquitous in cellular pathways 
(\cite{Gardner, pomerening-ferrell,Selkov, Sha}).  In either case, the system under
meaningful biological parameters shows convergence, not other dynamical
properties such as periodic behavior or even chaotic behavior.  Analytical
studies done for the quasi-steady-state version of the model (slow dynamics),
which is a monotone system, indicate that the reduced system is indeed
monostable or bistable, see \cite{Kholodenko_steady_state}.
Thus, it is of great interest to show that, at least in certain parameter
ranges (as required by singular perturbation theory), the full system
inherits convergence properties from the reduced system, and this is what we
do as an application of our results.

A feature of our approach, as in other control problems
 \cite{Abed}, \cite{Isidori},
is the use of geometric invariant manifold theory
\cite{Fenichel,Jones,Nipp}. There is a manifold $M_\varepsilon $, invariant
for the full dynamics of a singularly perturbed system, which attracts all
near-enough solutions.
However, we need to exploit the full power of the theory, and especially the
fibration structure and an asymptotic phase property. The system
restricted to the invariant manifold $M_\varepsilon $ is a regular
perturbation of the slow ($\varepsilon $$=$$0$) system.  As remarked in
Theorem 1.2 in Hirsch's early paper~\cite{Hirsch2}, a $C^1$ regular
perturbation of a flow with eventually positive derivatives also has generic
convergence properties. 
So, solutions in the manifold will generally be well-behaved, and
asymptotic phase implies that solutions near $M_\varepsilon$ track solutions in $M_\varepsilon $,
and hence also converge to equilibria if solutions on $M_\varepsilon $ do. A
key technical detail is to establish that the tracking solutions also start
from the ``good'' set of initial conditions, for generic solutions of the
large system.

A preliminary version of these results was presented at the
2006 Conference on Decision and Control, and dealt with the special case
of singularly perturbed systems of the form:
\begin{align*}
\dot{x}=&f(x,y)\\
\varepsilon \dot{y}=&Ay+h(x) 
\end{align*}
on a product domain, where $A$ is a constant Hurwitz matrix and the reduced system
$\dot{x}=f(x,-A^{-1}h(x))$ is strongly monotone. However, for the application to the above futile cycle, there are two major
problems with that formulation: first, the dynamics of the fast system have to
be allowed to be nonlinear in $y$, and second, it is crucial to allow for an
$\varepsilon$-dependence on the right-hand side as well as to allow the domain
to be a polytope depending on $\varepsilon$.  We provide a much more general formulation here.

We note that no assumptions are imposed regarding global convergence of the
reduced system, which is essential because of the intended application to
multi-stable systems.  This seems to rule out the applicability of
Lyapunov-theoretic and ISS tools~\cite{teel1,teel2}.

This paper is organized as follows. The main result is stated in Section
\ref{state}. In Section \ref{pre}, we review some basic definitions and
theorems about monotone systems. The detailed proof of the main theorem can be
found in Section \ref{proof}, and applications to the MAPK system and another
set of ordinary differential equations are discussed in Section
\ref{application}. Finally, in Section \ref{conclusion}, we summarize the key
points of this paper.

\section{Statement of the Main Theorem}
\label{state}
In this paper, we focus on the dynamics of the following prototypical system
in singularly perturbed form: 
\begin{align}
\label{eqn:slow_0}
\frac{dx}{dt}&= f_0(x,y,\varepsilon) \\
\varepsilon \frac{dy}{dt}&= g_0(x,y,\varepsilon).\notag
\end{align}
We will be interested in the dynamics of this system on a time-varying domain $D_\varepsilon$. For $0<\varepsilon\ll1$, the variable $x$ changes much
slower than $y$. 
As long as $\varepsilon\not=0$, one may also change the time scale to
$\tau=t/\varepsilon$, and study the equivalent form:
\begin{align}
\label{eqn:f_xy}
\frac{dx}{d\tau}&=\varepsilon f_0(x,y,\varepsilon) \\
\frac{dy}{d\tau}&= g_0(x,y,\varepsilon). \notag
\end{align}
Within this general framework, we will
make the following assumptions (some technical terms will be defined later), 
where the integer $r>1$ and the positive number $\varepsilon_0$ are fixed from now on:

\begin{description}
\item[{\bf A1}] Let $U \subset \mathbb R^n$ and $V\subset \mathbb R^m$ be open and bounded. The functions 
\[
f_0:U \times V \times [0, \varepsilon_0] \rightarrow \mathbb R^n
\]
\[
g_0:U \times V \times [0, \varepsilon_0] \rightarrow \mathbb R^m
\]
are both of class $C^r_b$, where a function $f$ is in  $C_b^r$
if it is in $C^r$ and its derivatives up to order $r$ as well as $f$ itself
are bounded. 

\item[{\bf A2}] There is a function 
\[
m_0:U \rightarrow V
\]
of class $C_b^r$, such that $g_0(x,m_0(x),0)=0$ for all $x$ in $U$.

\end{description}

It is often helpful to consider $z=y-m_0(x)$, and the fast system \eqref{eqn:f_xy} in the new coordinates becomes:
\begin{align}
\label{eqn:f_xz}
\frac{dx}{d\tau}&=\varepsilon f_1(x,z,\varepsilon)  \\
\frac{dz}{d\tau}&= g_1(x,z,\varepsilon), \notag
\end{align}
where
\begin{align*}
f_1(x,z,\varepsilon)&=f_0(x,z+m_0(x),\varepsilon),\\
g_1(x,z,\varepsilon)&=g_0(x,z+m_0(x),\varepsilon)-\varepsilon [D_x m_0(x)] f_1(x,z,\varepsilon).
\end{align*}
When $\varepsilon=0$, the system \eqref{eqn:f_xz} degenerates to
\begin{equation}
\label{eqn:f_z}
\frac{dz}{d\tau}=g_1(x,z,0), \ \ x(\tau) \equiv x_0 \in U,
\end{equation}
seen as equations on $\{z \,|\,z+m_0(x_0) \in V\}$.
\begin{description}
\item[{\bf A3}] 
The steady state $z=0$ of \eqref{eqn:f_z} is globally asymptotically stable on $\{z \,|\,z+m_0(x_0) \in V\}$ for all $x_0 \in U$. 
\end{description}
\begin{description}
\item[{\bf A4}] All eigenvalues of the matrix $D_y g_0(x,m_0(x),0)$ have negative real parts for every $x \in U$, i.e. the matrix $D_y g_0(x,m_0(x),0)$ is Hurwitz on $U$.

\item[{\bf A5}] There exists a family of convex compact sets $D_\varepsilon
  \subset U \times V$, which depend continuously on 
$\varepsilon \in[0,\varepsilon_0]$, such that \eqref{eqn:slow_0} is positively
  invariant on $D_\varepsilon$ for $\varepsilon \in (0,\varepsilon_0]$. 

\item[{\bf A6}] The flow $\psi_t^0$ of the limiting system (set $\varepsilon=0$ in \eqref{eqn:slow_0}):
\begin{align}
\label{eqn:m_0}
\frac{dx}{dt}= f_0(x,m_0(x),0)
\end{align}
has eventually positive derivatives on $K_0$, where $K_0$ is the projection of
\[
D_0 \bigcap \{(x,y)\,|\,y=m_0(x),x \in U\}
\]
onto the $x$-axis. 

\item[{\bf A7}] 
The set of equilibria of \eqref{eqn:slow_0} on $D_\varepsilon$ is totally disconnected. 
\end{description}
\begin{rk}
Assumption {\bf A3} implies that $y=m_0(x)$ is a unique solution of $g_0(x,y,0)=0$ on $U$.

Continuity in {\bf A5} is understood with respect to the Hausdorff metric. 

In mass-action chemical kinetics, the vector fields are polynomials. 
So, {\bf A1} follows naturally. 
\end{rk}

Our main theorem is:
\begin{thm}
\label{thm:main}
Under assumptions {\bf A1} to {\bf A7}, there exists a positive constant $\varepsilon^*< \varepsilon_0$ such that for each $\varepsilon \in (0,\varepsilon^*)$, the forward trajectory of \eqref{eqn:slow_0} starting from almost every point in $D_\varepsilon$ converges to some equilibrium. 
\end{thm}

\section{Monotone Systems of Ordinary Differential Equations}
\label{pre}

In this section, we review several useful definitions and theorems regarding
monotone systems.  As we wish to provide results valid for arbitrary orders,
not merely orthants, and some of these results, though well-known, are not
readily available in a form needed for reference, we provide some technical
proofs.

\begin{define}
A nonempty, closed set $C \subset \mathbb R^N$ is a cone if 
\begin{enumerate}
\item $C+C \subset C$,
\item $\mathbb R_+ C \subset C$,
\item $C \bigcap \,(-C)= \{0\}$.
\end{enumerate}
\end{define}
We always assume $C \not= \{0\}$. 
Associated to a cone $C$ is a partial order on $\mathbb R^N$. For any $x,y
\in \mathbb R^N$, we define:
\begin{align*}
x \geq y & \Leftrightarrow x-y \in C \\
x > y & \Leftrightarrow x-y \in C, x \not= y.
\end{align*}
When Int$C$ is not empty, we can define
\[
x \gg y \Leftrightarrow x-y \in \mbox{Int}C.
\]
\begin{define}
The dual cone of $C$ is defined as
\[
C^*=\{\lambda \in (\mathbb R^N)^*\, |\, \lambda(C) \geq 0\}.
\]
\end{define}
An immediate consequence is
\begin{align*}
x \in C & \Leftrightarrow \lambda(x) \geq 0, \forall \lambda \in C^* \\
x \in \mbox{Int}C & \Leftrightarrow \lambda(x)>0, \forall \lambda \in C^* \setminus \{0\}.
\end{align*}

With this partial ordering on $\mathbb R^N$, we analyze certain features of
the dynamics of an ordinary differential equation: 
\begin{eqnarray}
\frac{d z}{dt}=F(z), \label{eqn:def}
\end{eqnarray}
where $F:\mathbb R^N \rightarrow \mathbb R^N$ is a $C^1$ vector field. We are interested in a special class of equations which preserve the ordering along the trajectories. For simplicity, the solutions of \eqref{eqn:def} are assumed to exist for all $t \geq 0$ in the sets considered below.

\begin{define}
The flow $\phi_t$ of (\ref{eqn:def}) is said to have (eventually) positive derivatives on a set $V \subseteq \mathbb R^N$, if $[D_z \phi_t( z)]x \in \mbox{Int}C$ for all $x \in C \setminus \{0\},z \in V$, and $t  \geq 0$ ($t \geq t_0$ for some $t_0>0$). 
\end{define}

It is worth noticing that $[D_z \phi_t(z)]x \in \mbox{Int}C$ is equivalent to $\lambda([D_z \phi_t( z)]x)>0$ for all $\lambda \in C^*$ with norm one. We will use this fact in the proof of the next lemma, which deals with ``regular'' perturbations in the dynamics. The proof is in the same spirit as in Theorem 1.2 of \cite{Hirsch}, but generalized to the arbitrary cone $C$.

\begin{lemma}
Assume $V \subset \mathbb R^N$ is a compact set in which the flow $\phi_t$ of \eqref{eqn:def} has eventually positive derivatives. Then there exists $\delta>0$ with the following property. Let $\psi_t$ denote the flow of a $C^1$ vector field $G$ such that the $C^1$ norm of $F(z)-G(z)$ is less than $\delta$ for all $z$ in $V$. Then there exists $t_*>0$ such that if $\psi_s(z) \in V$ for all $s \in [0,t]$ where $t\geq t_*$, then $[D_z \psi_t(z)]x \in \mbox{Int}C$ for all $z \in V$ and $x \in C \setminus \{0\}$.
\label{lemma:Hirsch}
\end{lemma}
\begin{proof}
Pick $t_0>0$ so that $\lambda([D_z \phi_t(z)]x)>0$ for all $t \geq t_0, z \in V, \lambda \in C^*, x \in C$ with $|\lambda|=1, |x|=1$. Then there exists $\delta>0$ with the property that when the $C^1$ norm of $F(z)-G(z)$ is less than $\delta$, we have $\lambda([D_z \psi_t(z)]x)>0$ for $t_0 \leq t \leq 2t_0$.

When $t>2t_0$, we write $t=r+kt_0$, where $t_0 \leq r < 2t_0$ and $k \in \mathbb N$. If $\psi_s(z) \in V$ for all $s \in [0,t]$, we can define $z_j:=\psi_{jt_0}(z)$ for $j=0,\dots, k$. For any $x \in C \setminus \{0\}$, using the chain rule, we have:
\[
[D_z \psi_t(z)]x=[D_z \psi_r(z_k)][D_z \psi_{t_0}(z_{k-1})] \cdots [D_z \psi_{t_0}(z_0)]x.
\]
By induction, it is easy to see that $[D_z \psi_t(z)]x \in \mbox{Int}C$.
\end{proof}

\begin{cor}
\label{cor:Hirsch}
If $V$ is positively invariant under the flow $\psi_t$, then $\psi_t$ has eventually positive derivatives in $V$.
\end{cor}
\begin{proof} 
If $V$ is positively invariant under the flow $\psi_t$, then for any $z \in V$ the condition $\psi_s(z) \in V$ for $s \in [0,t]$ is satisfied for all $t \geq 0$. By the previous lemma, $\psi_t$ has eventually positive derivatives in $V$.
\end{proof}

\begin{define}
The system \eqref{eqn:def} or the flow $\phi_t$ of \eqref{eqn:def} is called monotone (resp. strongly monotone) in a set $W \subseteq \mathbb R^N$, if for all $t \geq 0$ and $z_1,z_2 \in W$,
\begin{align*}
z_1 \geq z_2 \Rightarrow &\phi_t(z_1) \geq \phi_t(z_2) \\ 
&(\mbox{resp.} \  \phi_t(z_1) \gg \phi_t(z_2) \ \mbox{when} \ z_1 \not= z_2).
\end{align*}
It is eventually (strongly) monotone if there exists $t_0 \geq 0$ such that $\phi_t$ is (strongly) monotone for all $t \geq t_0$.
\end{define}
\begin{define}
An set $W \subseteq \mathbb R^N$ is called p-convex, if $W$ contains the entire line segment joining $x$ and $y$ whenever $x \leq y$, $x,y \in W$.
\end{define}
The next two propositions discuss the relations between the two definitions, (eventually) positive derivatives and (eventually) strongly monotone.
\begin{pro}
\label{pro:EPD_ESM}
Let $W \subseteq \mathbb R^N$ be p-convex. If the flow $\phi_t$ has (eventually) positive derivatives in $W$, then it is (eventually) strongly monotone in $W$.
\end{pro}
\begin{proof}
For any $z_1>z_2 \in W, \lambda \in C^*\setminus \{0\}$ and $t \geq 0$ ($t
\geq t_0$ for some $t_0 >0$), we have that $\lambda(\phi_t(z_1)-\phi_t(z_2))$
equals
\[
\int_0^1 \lambda \big([D_z \phi_t(sz_1+(1-s)z_2)](z_1-z_2)\big) ds>0.
\]
Therefore, $\phi_t$ is (eventually) strongly monotone in $W$.
\end{proof}

\begin{pro}
\label{pro:ESM_EPD}
Suppose $\phi_t$ is (eventually) strongly monotone on an open set $U \subseteq \mathbb R^N$. Then $\phi_t$ has (eventually) positive derivatives in $U$.
\end{pro}
\begin{proof}
Fix $t>0$ such that $\phi_t$ as a function from $U$ to $\mathbb R^N$ is strongly monotone (i.e. $\phi_t(z_1) \gg \phi_t(z_2)$, whenever $z_1 \geq z_2$, $z_1 \not= z_2$ in $U$). 

For any $\lambda \in C^* \setminus \{0\},x \in C \setminus \{0\}$, and $z \in U$,
\[
\lambda([D_z\phi_t(z)]x)=\lim_{h \rightarrow 0}h^{-1}\lambda(\phi_t(z+hx)-\phi_t(z)).
\]
If $h>0$ is small enough such that $z+hx \in U$, then we have $\lambda \big(\phi_t(z+hx)- \phi_t(z)\big)>0$. Thus $\lambda([D_z\phi_t(z)]x)>0$, and $[D_z\phi_t(z)]x \in \mbox{Int}C$.
\end{proof}
\begin{lemma}
\label{lemma:Hirsch_conv}
Suppose that the flow $\phi_t$ of \eqref{eqn:def} has compact closure and eventually positive derivatives in a $p$-convex set $W \subseteq \mathbb R^N$. If the set of equilibria is totally disconnected (e.g. countable), then the forward trajectory starting from almost every point in $W$ converges to an equilibrium.
\end{lemma}
\begin{proof}
By Proposition \ref{pro:EPD_ESM}, $\phi_t$ is eventually strongly
monotone. The result easily follows from Hirsch's Generic
Convergence Theorem (\cite{Hirsch-Smith}, \cite{Smith}). 
\end{proof}

\section{Details of the Proof}
\label{proof}
Our approach to solve the varying domain problem is motivated by Nipp \cite{Nipp}.  The idea is to extend the vector field from $U \times V
\times [0,\varepsilon]$ to $\mathbb R^n \times \mathbb R^m \times
       [0,\varepsilon_0]$, then apply geometric singular perturbation theorems
       (\cite{Sakamoto}) on $\mathbb R^n \times \mathbb R^m \times
       [0,\varepsilon_0]$, and finally restrict the flows to $D_\varepsilon$ for
       the generic convergence result. 
\subsection{Extensions of the vector field}
\label{extension}
For a given compact set $K \subset \mathbb R^n$ ($K_0 \subseteq K \subset U$), the following procedure is adopted from \cite{Nipp} to extend a $C_b^r$ function with respect to the $x$ coordinate from $U$ to $\mathbb R^n$, such that the extended function is $C^r_b$ and agrees with the old one on $K$.
This is a routine ``smooth patching'' argument.

Let $U_1$ be an open subset of $U$ with $C^r$ boundary and such that $K \subset U_1 \subseteq U$. For $\Theta_0>0$ sufficiently small, define
\[
U_1^{\Theta_0}:=\{x \in U_1 \,|\, \Theta(x)\geq \Theta_0\}, \ \mbox{where} \ \Theta(x):=\min_{u \in \partial U_1} |x-u|,
\]
such that $K$ is contained in $U_1^{\Theta_0}$. Consider the scalar $C^\infty$ function $\rho$:
\[
\rho(a):=\left \{ 
\begin{array}{ll}
0& a \leq 0 \\
\exp(1-\exp(a-1)/a)& 0 <a<1 \\
1& a \geq 1.
\end{array} \right. 
\]
Define
\[
\hat \Theta(x):=\left \{ 
\begin{array}{ll}
0& x \in \mathbb R^n \setminus U_1 \\
\Theta(x) & x \in U_1 \setminus U_1^{\Theta_0}\\
\Theta_0 & x \in U_1^{\Theta_0},
\end{array} \right.
\]
and
\[
\bar \Theta (x):=\rho\big(\frac{\hat \Theta(x)}{\Theta_0}\big).
\]
For any $q \in C^r_b(U)$, let
\[
\bar{\bar q}(x):=\left \{
\begin{array}{ll}
q(x)& x \in U_1 \\
0& x \in \mathbb R^n \setminus U_1,
\end{array} \right.
\mbox{and} \  \bar q(x):=\bar \Theta(x)\bar{\bar q}(x).
\]
Then $\bar q(x) \in C^r_b(\mathbb R^n)$ and $\bar q(x) \equiv q(x)$ on $K$.

We fix some $d_0>0$ such that
\[
D_{d_0}:=\{ z \in \mathbb R^m \, | \, |z| \leq d_0\} \subset \bigcap_{x \in K} \{z \,|\, z+m_0(x) \in V\}.
\]
Then we extend the functions $f_1$ and $m_0$ to $\bar f_1$ and $\bar m_0$ respectively with respect to $x$ in the above way. To extend $g_1$, let us first rewrite the differential equation for $z$ as:
\[
\frac{dz}{d\tau}\!=\![B(x)+C(x,z)]z+\varepsilon H(x,z,\varepsilon)-\varepsilon [D_x m_0(x)]f_1(x,z,\varepsilon),
\]
where 
\[
B(x)=D_y g_0(x,m_0(x),0) \ \mbox{and} \ C(x,0)=0.
\]
Following the above procedures, we extend the functions $C$ and $H$ to $\bar C$ and $\bar H$, but the extension of $B$ is defined as
\[
\bar B(x):=\bar \Theta(x)\bar {\bar B}(x)-\mu(1-\bar \Theta(x))I_n,
\]
where $\mu$ is the positive constant such that the real parts of all eigenvalues of $B(x)$ is less than $-\mu$ for every $x \in K$. According to the definition of $\bar B(x)$, all eigenvalues of $\bar B(x)$ will have negative real parts less than $-\mu$ {\em for every} $x \in \mathbb R^n$. The extension $\bar g_1$, defined as:
\[
[\bar B(x)+\bar C(x,z)]z+\varepsilon \bar H(x,z,\varepsilon)-\varepsilon [D_x \bar m_0(x)] \bar f_1(x,z,\varepsilon),
\]
is then $C^{r-1}_b(\mathbb R^n \times D_{d_0} \times [0,\varepsilon_0])$ and agrees with $g_1$ on $K \times D_{d_0} \times [0, \varepsilon_0]$.

To extend functions $\bar f_1$ and $\bar g_1$ in the $z$ direction from $D_{d_0}$ to $\mathbb R^m$, we use the same extension technique but with respect to $z$. Let us denote the extensions of $\bar f_1, \bar C, \bar H$ and the function $z=z$ by $\tilde f_1, \tilde C, \tilde H$ and $\tilde z$ respectively, then define
$\tilde g_1$ as:
\[
[\bar B(x)+\tilde C(x,z)]\tilde z(z)+\varepsilon \tilde H(x,z,\varepsilon)-\varepsilon [D_x \bar m_0(x)] \tilde f_1(x,z,\varepsilon),
\]
which is now $C^{r-1}_b(\mathbb R^n \times \mathbb R^m \times [0,\varepsilon_0])$ and agrees with $g_1$ on $K \times D_{d_1} \times [0, \varepsilon_0]$ for some $d_1$ slightly less than $d_0$. Notice that $z=0$ is a solution of $\tilde g_1(x,z,0)=0$, which guarantees that for the extended system in $(x,y)$ coordinates ($y=z+\bar m_0(x)$):
\begin{align}
\label{eqn:f_e_xy}
\frac{dx}{d\tau}&=\varepsilon f(x,y,\varepsilon) \\
\frac{dy}{d\tau}&=g(x,y,\varepsilon), \notag
\end{align} 
$y=\bar m_0(x)$ is the solution of $g(x,y,0)=0$.
To summarize, \eqref{eqn:f_e_xy} satisfies
\begin{description}
\item[{\bf E1}] The functions 
\[
f \in C^r_b(\mathbb R^n \times \mathbb R^m \times [0,\varepsilon_0]), 
\]
\[
g \in C^{r-1}_b(\mathbb R^n \times \mathbb R^m \times [0,\varepsilon_0]), 
\]
\[
\bar m_0 \in C^r_b(\mathbb R^n),\ \  g(x,\bar m_0(x),0)=0, \ \forall x \in \mathbb R^n.
\]

\item[{\bf E2}] All eigenvalues of the matrix $D_y g(x,\bar m_0(x),0)$ have negative real parts less than $-\mu$ for every $x \in \mathbb R^n$.
\item[{\bf E3}] The function $\bar m_0$ coincides with $m_0$ on $K$, and the functions $f$ and $g$ coincide with $f_0$ and $g_0$ respectively on 
\[
\Omega_{d_1}:=\{(x,y) \, | \, x \in K, \, |y-m_0(x)| \leq d_1 \}.
\]

\end{description}
Conditions {\bf E1} and {\bf E2} are the assumptions for geometric singular
perturbation theorems, and condition {\bf E3} ensures that on $\Omega_{d_1}$
the flow of \eqref{eqn:f_xy} coincides with the flow of \eqref{eqn:f_e_xy}. If
we apply geometric singular perturbation theorems to \eqref{eqn:f_e_xy} on
$\mathbb R^n \times \mathbb R^m \times [0,\varepsilon_0]$, the exact same
results are true for \eqref{eqn:f_xy} on $\Omega_{d_1}$. For the rest of the
paper, we identify the flow of \eqref{eqn:f_e_xy} and the flow of
\eqref{eqn:f_xy} on $\Omega_{d_1}$ without further mentioning this fact.(Later, in Lemmas \ref{lemma:positive}-\ref{l:intersect}, when globalizing the results, we consider again the original system.) 

\subsection{Geometric singular perturbation theory}

The theory of geometric singular perturbation can be traced back to the work
of Fenichel \cite{Fenichel}, which first revealed the geometric aspects
of singular perturbation problems. Later on, the works by Knobloch and Aulbach
\cite{Knobloch_Aulbach}, Nipp \cite{Nipp}, and Sakamoto \cite{Sakamoto} also
presented results similar to \cite{Fenichel}. By now, the theory is fairly
standard, and there have been enormous applications to traveling waves of
partial differential equations, see \cite{Jones} and the references there. For
control theoretic applications, see \cite{Abed,Isidori}. 

To apply geometric singular perturbation theorems to the vector field on
$\mathbb R^n \times \mathbb R^m \times[0,\varepsilon_0]$, we use the theorems
stated in \cite{Sakamoto}. 
The following lemma is a restatement of the theorems in \cite{Sakamoto}, and we refer to \cite{Sakamoto} for the proof.

\begin{lemma}
\label{l:Fenichel}
Under conditions {\bf E1} and {\bf E2}, there exists a positive $\varepsilon_1<\varepsilon_0$ such that for every $\varepsilon \in (0,\varepsilon_1]$:
\begin{enumerate}
\item There is a $C^{r-1}_b$ function
\[
m: \, \mathbb R^n \times [0,\varepsilon_1] \rightarrow \mathbb R^m
\]
such that the set $M_\varepsilon$ defined by 
\[
M_\varepsilon:=\{\big(x,m(x,\varepsilon)\big)\,|\, x \in \mathbb R^n\}
\]
is invariant under the flow generated by (\ref{eqn:f_e_xy}). Moreover,
\[
\sup_{x \in \mathbb R^n} |m(x,\varepsilon)-\bar m_0(x)|=O(\varepsilon), \ \mbox{as}\ \varepsilon \rightarrow 0.
\]
In particular, we have $m(x,0)=\bar m_0(x)$ for all $x \in \mathbb R^n$.
\item The set consisting of all the points $(x_0, y_0)$ such that 
\[
\sup_{\tau \geq 0} |y(\tau; x_0,y_0)-m\big(x(\tau; x_0,y_0),\varepsilon\big)|e^{\frac{\mu\tau}{4}}<\infty,
\]
where $\big(x(\tau; x_0,y_0), y(\tau; x_0, y_0)\big)$ is the solution of \eqref{eqn:f_e_xy} passing through $(x_0,y_0)$ at $\tau=0$, is a $C^{r-1}$-immersed submanifold in $\mathbb R^n \times \mathbb R^m$ of dimension $n+m$, denoted by $W^s(M_\varepsilon)$.

\item There is a positive constant $\delta_0$ such that if
\[
\sup_{\tau \geq 0} |y(\tau; x_0,y_0)-m\big(x(\tau;x_0,y_0),\varepsilon\big)|<\delta_0,
\]
then $(x_0,y_0) \in W^s(M_\varepsilon)$.

\item The manifold $W^s(M_\varepsilon)$ is a disjoint union of $C^{r-1}$-immersed manifolds $W^s_\varepsilon(\xi)$ of dimension $m$:
\[
W^s(M_\varepsilon)=\bigcup_{\xi \in \mathbb R^n} W_\varepsilon^s(\xi).
\]
For each $\xi \in \mathbb R^n$, let $H_\varepsilon(\xi)(\tau)$ be the solution for $\tau \geq 0$ of
\[
\frac{dx}{d\tau}= \varepsilon f(x,m(x,\varepsilon),\varepsilon), \ \ x(0)=\xi \in \mathbb R^n.
\]
Then, the manifold $W^s_\varepsilon(\xi)$ is the set
\[
\{(x_0,y_0)\, | \, \sup_{\tau \geq 0} |\tilde x(\tau)|e^{\frac{\mu\tau}{4}}<\infty,\sup_{\tau \geq 0} |\tilde y(\tau)|e^{\frac{\mu\tau}{4}}<\infty \},
\]
where 
\begin{align*}
\tilde x(\tau)&=x(\tau; x_0,y_0)-H_\varepsilon(\xi)(\tau), \\
\tilde y(\tau)&=y(\tau; x_0,y_0)-m\big(H_\varepsilon(\xi)(\tau),\varepsilon\big).
\end{align*}
\item The fibers are ``positively invariant'' in the sense that $W^s_\varepsilon(H_\varepsilon(\xi)(\tau))$ is the set
\[
\{\big(x(\tau;x_0,y_0), y(\tau;x_0,y_0)\big) \, | \, (x_0,y_0) \in W^s_\varepsilon(\xi)\}
\]
for each $\tau \geq 0$, see Figure \ref{fig:We_fiber}.
\item The fibers restricted to the $\delta_0$ neighborhood of $M_\varepsilon$,
  denoted by $W^s_{\varepsilon,\delta_0}$, can be parametrized as
  follows. There   are two $C^{r-1}_b$ functions 
\[
P_{\varepsilon,\delta_0}: \mathbb R^n \times D_{\delta_0} \rightarrow \mathbb R^n
\]
\[
Q_{\varepsilon,\delta_0}: \mathbb R^m \times D_{\delta_0} \rightarrow \mathbb R^m,
\]
and a map
\[
T_{\varepsilon,\delta_0}:\mathbb R^n \times D_{\delta_0} \rightarrow \mathbb R^n \times \mathbb R^m 
\]
mapping $(\xi,\eta)$ to $(x,y)$, where
\[
x=\xi+P_{\varepsilon,\delta_0}(\xi,\eta), \ \ y=m(x,\varepsilon)+ Q_{\varepsilon,\delta_0}(\xi,\eta)
\]
such that
\[
W^s_{\varepsilon, \delta_0}(\xi)=T_{\varepsilon,\delta_0}(\xi,D_{\delta_0}).
\]
\end{enumerate}
\end{lemma}
\begin{flushright}
\QED
\end{flushright}
\begin{rk}
The $\delta_0$ in property 3 can be chosen uniformly for $\varepsilon \in (0,\varepsilon_0]$. Without loss of generality, we assume that $\delta_0<d_1$.

Notice that property 4 insures that for each $(x_0,y_0)\in W^s(M_\varepsilon)$, there exists a $\xi$ such that 
\[
|x(\tau;x_0,y_0)-H_\varepsilon(\xi)(\tau)| \rightarrow 0,
\]
\[
|y(\tau; x_0,y_0)-m\big(H_\varepsilon(\xi)(\tau),\varepsilon\big)| \rightarrow 0.
\]
as $t \rightarrow 0$. This is often referred as the ``asymptotic phase'' property, see Figure \ref{fig:We_fiber}.
\begin{figure}[h]
  \centering \includegraphics[scale=0.5]{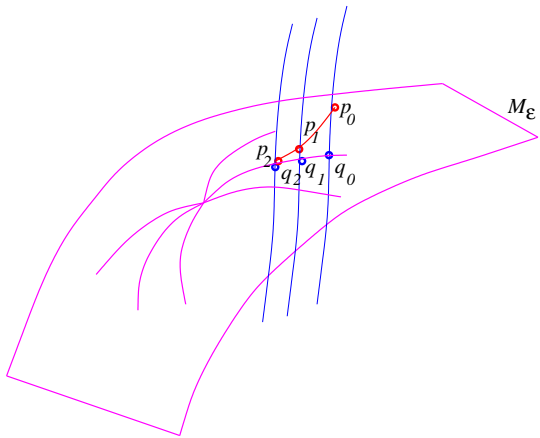}
  \caption{An illustration of the ``positive invariant'' and ``asymptotic phase'' properties. Let $p_0$ be a point on the fiber $W^s_\varepsilon(q_0)$ (vertical curve). Suppose the solution of \eqref{eqn:f_e_xy} starting from $q_0 \in M_\varepsilon$ evolves to $q_1 \in M_\varepsilon$ at time $\tau_1$, then the solution of \eqref{eqn:f_e_xy} starting from $p_0$ will evolve to $p_1 \in W^s_\varepsilon(q_1)$ at time $\tau_1$. At time $\tau_2$, they evolve to $q_2,p_2$ respectively. These two solutions are always on the same fiber. If we know that the one starting from $q_0$ converges to a equilibrium, then the one starting from $p_0$ also converges to a equilibrium.}
\label{fig:We_fiber}
\end{figure}
\end{rk}
\subsection{Further analysis of the dynamics}
The first property of Lemma \ref{l:Fenichel} concludes the existence of an invariant manifold $M_\varepsilon$. 
There are two reasons to introduce $M_\varepsilon$. First, on $M_\varepsilon$ the $x$-equation is decoupled from the $y$-equation:
\begin{align}
\label{eqn:m_e}
\frac{dx}{dt}&= f(x,m(x,\varepsilon),\varepsilon) \\
y(t)&= m(x(t), \varepsilon). \notag
\end{align}
This reduction allows us to analyze a lower dimensional system, whose dynamics
may have been well studied. Second, when $\varepsilon$ approaches zero, the
limit of \eqref{eqn:m_e} is \eqref{eqn:m_0} on $K_0$. If \eqref{eqn:m_0} has
some desirable property, it is natural to expect that this property is
inherited by \eqref{eqn:m_e}. An example of this principle is provided by the
following Lemma: 
\begin{lemma}
\label{lemma:positive}
There exists a positive constant $\varepsilon_2<\varepsilon_1$, such that for each $\varepsilon \in (0,\varepsilon_2)$, the flow $\psi^\varepsilon_t$ of \eqref{eqn:m_e} has eventually positive derivatives on $K_\varepsilon$, which is the projection of $M_\varepsilon \bigcap D_\varepsilon$ to 
the $x$-axis.
\end{lemma}
\begin{proof}
Assumption {\bf A6} states that the flow $\psi_t^0$ of the limiting system \eqref{eqn:m_0} has eventually positive derivatives on $K_0$. By the continuity of $m(x,\varepsilon)$ and $D_\varepsilon$ at $\varepsilon$$=$$0$, we can pick $\varepsilon_2$ small enough such that the flow $\psi_t^0$ has eventually positive derivatives on $K_\varepsilon$ for all $\varepsilon \in (0,\varepsilon_2)$. Applying Corollary \ref{cor:Hirsch}, we conclude that the flow 
$\psi^\varepsilon_t$ of \eqref{eqn:m_e} has eventually positive derivatives on $K_\varepsilon$ provided $K_\varepsilon$ is positively invariant under (\ref{eqn:m_e}), which follows easily from the facts that (\ref{eqn:f_e_xy}) is positively invariant on $D_\varepsilon$ and $M_\varepsilon$ is an invariant manifold.
\end{proof}

\begin{figure}[h]
  \centering \includegraphics[scale=0.4]{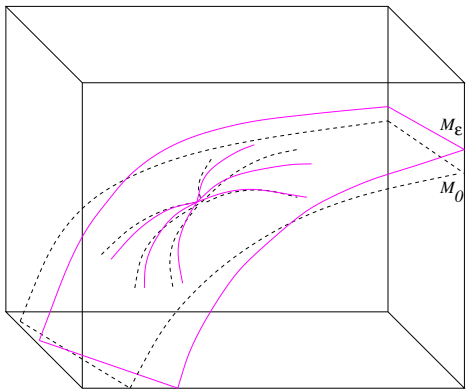}
  \caption{This is a sketch of the manifolds $M_0$ (surface bounded by dashed curves), $M_\varepsilon$ (surface bounded by dotted curves), and $D_\varepsilon$ (the cube). It highlights two major characters of $M_\varepsilon$. First, $M_\varepsilon$ is close to $M_0$. Second, the trajectories on $M_\varepsilon$ converge to 
equilibria if those on $M_0$ do.}
\label{fig:M_e}
\end{figure}

The next Lemma asserts that the generic convergence property is preserved for \eqref{eqn:m_e}, see Figure \ref{fig:M_e}.

\begin{lemma}
\label{l:conv}
For each $\varepsilon \in (0,\varepsilon_2)$, there exists a set $C_\varepsilon \subseteq K_\varepsilon$ such that the forward trajectory of \eqref{eqn:m_e} starting from any point of $C_\varepsilon$ converges to some equilibrium, and the Lebesgue measure of $K_\varepsilon \setminus C_\varepsilon$ is zero.
\end{lemma}
\begin{proof}
Apply Lemma \ref{lemma:Hirsch_conv} and Lemma \ref{lemma:positive} under assumptions {\bf A5} and {\bf A7}.
\end{proof}

By now, we have discussed flows restricted to the invariant manifold $M_\varepsilon$. Next, we will explore the conditions for a point to be on $W^s(M_\varepsilon)$, the stable manifold of $M_\varepsilon$. Property 3 of Lemma \ref{l:Fenichel} provides a sufficient condition, namely, any point $(x_0,y_0)$ such that 
\begin{equation}
\label{eqn:Ws}
\sup_{\tau \geq 0} |y(\tau; x_0,y_0)-m\big(x(\tau;x_0,y_0),\varepsilon\big)|<\delta_0
\end{equation}
is on $W^s(M_\varepsilon)$. In fact, if we know that the difference between
$y_0$ and $m(x_0,\varepsilon)$ is sufficiently small, then the above condition
is always satisfied. More precisely, we have:

\begin{lemma}
\label{l:cond}
There exists $\varepsilon_3>0,\delta_0>d>0$, such that for each $\varepsilon \in (0,\varepsilon_3)$, if the initial condition satisfies $|y_0-m(x_0,\varepsilon)|<d$, then \eqref{eqn:Ws} holds, i.e. $(x_0,y_0) \in W^s(M_\varepsilon)$. 
\end{lemma}
\begin{proof}
Follows from the proof of Claim 1 in \cite{Nipp}.
\end{proof}

Before we get further into the technical details, let us give an outline of the proof of the main theorem. The proof can be decomposed into three steps. First, we show that almost every trajectory on $D_\varepsilon \bigcap M_\varepsilon$ converges to some equilibrium. This is precisely Lemma \ref{l:conv}. Second, we show that almost every trajectory starting from $W^s(M_\varepsilon)$ converges to some equilibrium. This follows from Lemma \ref{l:conv} and the ``asymptotic phase'' property in Lemma \ref{l:Fenichel}, but we still need to show that the set of non-convergent initial conditions is of measure zero. The last step is to show that all trajectories in $D_\varepsilon$ will eventually stay in $W^s(M_\varepsilon)$, which is our next lemma:
\begin{lemma}
\label{l:intersect}
There exist positive $\tau_0$ and $\varepsilon_4<\varepsilon_3$, such that $(x(\tau_0),y(\tau_0)) \in W^s(M_\varepsilon)$ for all $\varepsilon \in (0,\varepsilon_4)$, where $(x(\tau),y(\tau))$ is the solution to \eqref{eqn:f_xy} with the initial condition $(x_0,y_0) \in D_\varepsilon$.
\end{lemma}
\begin{proof}
It is convenient to consider the problem in $(x,z)$ coordinates. Let $(x(\tau),z(\tau))$ be the solution to \eqref{eqn:f_xz} with initial condition $(x_0,z_0)$, where $z_0=y_0-m(x_0,0)$. We first show that there exists a $\tau_0$ such that $|z(\tau_0)| \leq d/2$.

Expanding $g_1(x,z,\varepsilon)$ at the point $(x_0,z,0)$, the equation of $z$ becomes 
\[
\frac{dz}{d\tau} =g_1(x_0,z,0)+\frac{\partial g_1}{\partial x}(\xi,z,0)(x-x_0)+\varepsilon R(x,z,\varepsilon)
\]
for some $\xi(\tau)$ between $x_0$ and $x(\tau)$ (where $\xi(\tau)$ can be picked continuously in $\tau$). Let us write 
\[
z(\tau)=z^0(\tau)+w(\tau),
\]
where $z^0(\tau)$ is the solution to \eqref{eqn:f_z} with initial the condition $z^0(0)=z_0$, and $w(\tau)$ satisfies
\begin{align}
\label{eqn:w}
\frac{d w}{d\tau}&=g_1(x_0,z,0)-g_1(x_0,z^0,0)+\frac{\partial g_1}{\partial x}(\xi,z,0)(x-x_0)+\varepsilon R(x,z,\varepsilon) \\
&=\frac{\partial g_1}{\partial z}(x_0,\zeta,0)w +\varepsilon\frac{\partial g_1}{\partial x}(\xi,z,0)\int_0^\tau f_1(x(s),z(s),\varepsilon)\,ds +\varepsilon R(x,z,\varepsilon), \notag
\end{align}
with the initial condition $w(0)=0$ and some $\zeta(\tau)$ between $z^0(\tau)$ and $z(\tau)$ (where $\zeta(\tau)$ can be picked continuously in $\tau$). 

By assumption {\bf A3}, there exist a positive $\tau_0$ such that $|z^0(\tau)| \leq d/4$ for all $\tau \geq \tau_0$. Notice that we are working on the compact set $D_\varepsilon$, so $\tau_0$ can be chosen uniformly for all initial conditions in $D_\varepsilon$. 

We write the solution of \eqref{eqn:w} as:
\begin{align*}
w(\tau)&=\int_0^\tau\frac{\partial g_1}{\partial z}(x_0,\zeta,0)w \,ds+\varepsilon\int_0^\tau\!\!\big(\frac{\partial g_1}{\partial x}(\xi,z,0)\int_0^{s'}\!\!\! f_1(x,z,\varepsilon)\,ds'+R(x,z,\varepsilon)\big)\,ds.
\end{align*}
Since the functions $f_1, R$ and the derivatives of $g_1$ are bounded on $D_\varepsilon$, we have:
\[
|w(\tau)| \leq \int_0^\tau L|w| \,ds +
\varepsilon\int_0^\tau\left(M_1\int_0^{s'} M_2 \,ds'+M_3\right)\,ds,
\]
for some positive constants $L, M_i$, $i=1,2,3$. The notation $|w|$ means the Euclidean norm of $w \in \mathbb R^m$. Moreover, if we define
\[
\alpha(\tau)=\int_0^{\tau}\left(M_1\int_0^{s'} M_2 \,ds'+M_3\right)\,ds,
\]
then 
\[
|w(\tau)| \leq \int_0^\tau L|w|\,ds+\varepsilon \alpha(\tau_0), 
\]
for all $\tau \in [0,\tau_0]$ as $\alpha$ is increasing in $\tau$. Applying Gronwall's inequality (\cite{mct}), we have:
\[
|w(\tau)| \leq \varepsilon \alpha(\tau_0)e^{L\tau},
\]
which holds in particular at $\tau=\tau_0$. Finally, we choose $\varepsilon_4$ small enough such that $\varepsilon \alpha(\tau_0)e^{L\tau_0} < d/4$ and $|m(x,\varepsilon)-m(x,0)|<d/2$ for all $\varepsilon \in (0,\varepsilon_4)$. Then we have:
\begin{align*}
|y(\tau_0)-m(x(\tau_0),\varepsilon)| &\leq |y(\tau_0)-m(x(\tau_0),0)|+|m(x(\tau_0),\varepsilon)-m(x(\tau_0),0)|  \\
&< |z(\tau_0)|+d/2 \\
&<d/2+d/2=d.
\end{align*}
That is, $(x(\tau_0),y(\tau_0)) \in W^s(M_\varepsilon)$ by Lemma \ref{l:cond}. 
\end{proof}

By now, we have completed all these three steps, and are ready to prove Theorem \ref{thm:main}.

\subsection{Proof of Theorem \protect{\ref{thm:main}}}
\begin{proof}
Let $\varepsilon^*=\min\{\varepsilon_2,\varepsilon_4\}$. For $\varepsilon \in (0,\varepsilon^*)$, it is equivalent to prove the result for the fast system \eqref{eqn:f_xy}. Pick an arbitrary point $(x_0,y_0)$ in $D_\varepsilon$, and there are three cases:

\begin{enumerate}
\item $y_0=m(x_0,\varepsilon)$, that is, $(x_0,y_0)\in M_\varepsilon \bigcap D_\varepsilon$. By Lemma \ref{l:conv}, the forward trajectory converges to an equilibrium except for a set of measure zero.

\item $0<|y_0-m(x_0,\varepsilon)| < d$. By Lemma \ref{l:cond}, we know that $(x_0, y_0)$ is in $W^s(M_\varepsilon)$. Then, property 4 of Lemma \ref{l:Fenichel} guarantees that the point $(x_0,y_0)$ is on some fiber $W_{\varepsilon,d}^s(\xi)$, where $\xi \in K_\varepsilon$. If $\xi \in C_\varepsilon$, that is, the forward trajectory of $\xi$ converges to some equilibrium, then by the ``asymptotic phase'' property of Lemma \ref{l:Fenichel}, the forward trajectory of $(x_0,y_0)$ also converges to an equilibrium. To deal with the case when $\xi$ is not in $C_\varepsilon$, it is enough to show that the set 
\[
B_{\varepsilon,d}=\bigcup_{\xi \in K_\varepsilon \setminus C_\varepsilon} W^s_{\varepsilon,d}(\xi)
\] 
has measure zero in $\mathbb R^{m+n}$. Define 
\[
S_{\varepsilon,d}=(K_\varepsilon \setminus C_\varepsilon) \times D_{d}.
\]
By Lemma \ref{l:conv}, $K_\varepsilon \setminus C_\varepsilon$ has measure zero in $\mathbb R^n$, thus $S_{\varepsilon,d}$ has measure zero in $\mathbb R^n \times \mathbb R^m$. On the other hand, property 6 in Lemma \ref{l:Fenichel} implies $B_{\varepsilon,d}=T_{\varepsilon,d}(S_{\varepsilon,d})$. Since Lipschitz maps send measure zero sets to measure zero sets, $B_{\varepsilon,d}$ is of measure zero.
\item $|y_0-m(x_0,\varepsilon)| \geq d$. By Lemma \ref{l:intersect}, the point $\big(x(\tau_0),y(\tau_0)\big)$ is in $W^s(M_\varepsilon)$ and we are back to case 2. The proof is completed if the set $\phi^\varepsilon_{-\tau_0}(B_{\varepsilon,d})$ has measure zero, where $\phi^\varepsilon_\tau$ is the flow of \eqref{eqn:f_xy}. This is true because $\phi^\varepsilon_\tau$ is a diffeomorphism for any finite $\tau$. 
\end{enumerate}
\end{proof}

\section{Applications}
\label{application}
\subsection{An application to the dual futile cycle}
We assume that the reactions in Figure \ref{fig:Scheme} follow the
usual enzymatic mechanism (\cite{Ferrell_Bhatt}): 
\[
S_0  + E  \arrowschem{k_1}{k_{-1}}  C_1 \stackrel{k_2}{\rightarrow} S_1+ E \arrowschem{k_3}{k_{-3}} C_2 \stackrel{k_4}{\rightarrow} S_2+E
\]
\[
S_2+F \arrowschem{h_1}{h_{-1}} C_3 \stackrel{h_2}{\rightarrow} S_1+F \arrowschem{h_3}{h_{-3}} C_4 \stackrel{h_4}{\rightarrow} S_0+F.
\]
There are three conservation relations:
\begin{align*}
S_{tot}&=[S_0]+[S_1]+[S_2]+[C_1]+[C_2]+[C_4]+[C_3], \\
E_{tot}&=[E]+[C_1]+[C_2], \\
F_{tot}&=[F]+[C_4]+[C_3],
\end{align*}
where brackets indicate concentrations.
Based on mass action kinetics, we have the following set of ordinary differential equations:
\begin{align}
\label{eqn:ex_original}
\frac{d[S_0]}{d\tau}&=h_4[C_4] -k_1[S_0][E]+k_{-1}[C_1] \notag \\
 \frac{d[S_2]}{d\tau}&=k_4[C_2] -h_1[S_2][F]+h_{-1}[C_3]\notag \\
\frac{d[C_1]}{d\tau}&=k_1[S_0][E] -(k_{-1}+k_2)[C_1] \\
\frac{d[C_2]}{d\tau}&=k_3[S_1][E]-(k_{-3}+k_4)[C_2] \notag \\
\frac{d[C_4]}{d\tau}&=h_3[S_1][F]-(h_{-3}+h_4)[C_4] \notag \\
\frac{d[C_3]}{d\tau}&=h_1[S_2][F]-(h_{-1}+h_2)[C_3].  \notag
\end{align}
After rescaling the concentrations and time, \eqref{eqn:ex_original} becomes:
\begin{align}
\label{eqn:ex_full}
\frac{dx_1}{dt}&=-k_1S_{tot}x_1(1-y_1-y_2)+k_{-1}y_1+h_4 c y_3 \notag \\
\frac{dx_2}{dt}&=-h_1S_{tot}c x_2(1-y_3-y_4)+h_{-1}cy_4+k_4 y_2 \notag \\
\varepsilon\frac{dy_1}{dt}&=k_1S_{tot}x_1(1-y_1-y_2)-(k_{-1}+k_2)y_1  \\
\varepsilon\frac{dy_2}{dt}&=k_3S_{tot}(1-x_1-x_2-\varepsilon y_1-\varepsilon y_2-\varepsilon c y_3-\varepsilon c y_4) \times (1-y_1-y_2)  -(k_{-3}+k_4)y_2 \notag\\
\varepsilon\frac{dy_3}{dt}&=h_3S_{tot}(1-x_1-x_2-\varepsilon y_1-\varepsilon y_2-\varepsilon c y_3-\varepsilon c y_4) \times (1-y_3-y_4)-(h_{-3}+h_4)y_3 \notag\\
\varepsilon\frac{dy_4}{dt}&=h_1S_{tot}x_2(1-y_3-y_4)-(h_{-1}+h_2)y_4,\notag
\end{align}
where
\[
x_1=\frac{[S_0]}{S_{tot}}, \ \ x_2=\frac{[S_2]}{S_{tot}}, \ \ y_1=\frac{[C_1]}{E_{tot}}, \ \ y_2=\frac{[C_2]}{E_{tot}},
\]
\[
y_3=\frac{[C_4]}{F_{tot}}, \ \ y_4=\frac{[C_3]}{F_{tot}}, \ \ \varepsilon=\frac{E_{tot}}{S_{tot}}, \ \ c=\frac{F_{tot}}{E_{tot}},\ \ t=\tau \varepsilon.
\]
These equations are in the form of \eqref{eqn:slow_0}. 
The conservation laws suggest taking $\varepsilon_0=1/(1+c)$
and  
\begin{align*}
D_\varepsilon=&\{(x_1,x_2,y_1,y_2,y_3,y_4)\, | \, 0\leq y_1+y_2 \leq 1, 0\leq y_3+y_4 \leq 1, \\
&x_1 \geq 0, x_2 \geq 0, 0\leq x_1+x_2+\varepsilon(y_1+y_2+cy_3+cy_4) \leq 1\}.
\end{align*}
For $\varepsilon \in (0, \varepsilon_0]$, taking the inner product of the normal of $\partial D_\varepsilon$ and the vector field, it is easy to check that \eqref{eqn:ex_full} is positively invariant on $D_\varepsilon$, so {\bf A5} holds. We want to emphasize that in this example the domain $D_\varepsilon$ is a convex polytope varying with $\varepsilon$.

It can be proved that on $D_\epsilon$ system (12) has at most a finite number
of steady states, and thus {\bf A7} holds. This is a consequence of a more general result, proved
using some of the ideas given in~\cite{jeremy}, concerning the number of
steady states of more general systems of phosphorylation/dephosphorylation
reactions, see~\cite{ws3}.

Solving $g_0(x,y,0)=0$, we get
\begin{align*}
y_1=&\frac{x_1}{\frac{K_{m1}}{S_{tot}}+\frac{K_{m1}(1-x_1-x_2)}{K_{m2}}+x_1},\\
y_2=&\frac{\frac{K_{m1}(1-x_1-x_2)}{K_{m2}}}{\frac{K_{m1}}{S_{tot}}+\frac{K_{m1}(1-x_1-x_2)}{K_{m2}}+x_1}, \\
y_3=&\frac{\frac{K_{m3}(1-x_1-x_2)}{K_{m4}}}{\frac{K_{m3}}{S_{tot}}+\frac{K_{m3}(1-x_1-x_2)}{K_{m4}}+x_2}, \\
y_4=&\frac{x_2}{\frac{K_{m3}}{S_{tot}}+\frac{K_{m3}(1-x_1-x_2)}{K_{m4}}+x_2},
\end{align*}
where $K_{m1}, K_{m2}, K_{m3}$ and $K_{m4}$ are the Michaelis-Menten constants defined as
\[
K_{m1}=\frac{k_{-1}+k_2}{k_1}, \ \ K_{m2}=\frac{k_{-3}+k_4}{k_3}, \ \ K_{m3}=\frac{h_{-1}+h_2}{h_1}, \ \ K_{m4}=\frac{h_{-3}+h_4}{h_3}.
\]

Now, we need to find a proper set $U \subset \mathbb R^2$ satisfying assumptions {\bf A1}-{\bf A4}. Suppose that $U$ has the form
\[
U=\{(x_1,x_2)\, |\, x_1 > -\sigma , \, x_2>-\sigma,\, x_1+x_2 < 1+\sigma\},
\]
for some positive $\sigma$, and $V$ is any bounded open set such that $D_\varepsilon$ is contained in $U \times V$, then {\bf A1} follows naturally. Moreover, if
\[
\sigma \leq \sigma_0:=\min\left\{\frac{K_{m1}K_{m2}}{S_{tot}(K_{m1}+K_{m2})},\, \frac{K_{m3}K_{m4}}{S_{tot}(K_{m3}+K_{m4})}\right\},
\]
{\bf A2} also holds. To check {\bf A4}, let us look at the matrix:
\[
B(x):=D_y g_0(x,m_0(x),0)=\left( \begin{array}{cc}
         B_1(x) & 0 \\
         0 & B_2(x)
 \end{array} \right),
\]
where 
\[
B_1(x)= \left( \begin{array}{cc}
          -k_1S_{tot}x_1-(k_{-1}+k_2)  & -k_1S_{tot}x_1 \\
          -k_3S_{tot}(1-x_1-x_2) & -k_3S_{tot}(1-x_1-x_2)-(k_{-3}+k_4)
 \end{array} \right),
\]
and 
\[
B_2(x)=\left( \begin{array}{cc}
          -h_3S_{tot}(1-x_1-x_2)-(h_{-3}+h_4) & -h_3S_{tot}(1-x_1-x_2)\\
          -h_1S_{tot}x_2  & -h_1S_{tot}x_2-(h_{-1}+h_2
          \end{array} \right).
\]
If both matrices $B_1$ and $B_2$ have negative traces and positive determinants, then {\bf A4} holds.

Let us consider $B_1$ first. The trace of $B_1$ is \[
-k_1S_{tot}x_1-(k_{-1}+k_2)-k_3S_{tot}(1-x_1-x_2)-(k_{-3}+k_4).
\]
It is negative provided that 
\[
\sigma \leq \frac{k_{-1}+k_2+k_{-3}+k_4}{S_{tot}(k_1+k_3)}.
\]
The determinant of $B_1$ is
\begin{align*}
k_1(k_{-3}+k_4)S_{tot}x_1+k_3(k_{-1}+k_2)S_{tot}(1-x_1-x_2) +(k_{-1}+k_2)(k_{-3}+k_4).
\end{align*}
It is positive if 
\[
\sigma \leq \frac{(k_{-1}+k_2)(k_{-3}+k_4)}{S_{tot}\big(k_1(k_{-3}+k_4)+k_3(k_{-1}+k_2)\big)}.
\]
The condition for $B_2$ can be derived similarly. To summarize, if we take 
\begin{align*}
\sigma=\min\bigg\{& \sigma_0, \, \frac{k_{-1}+k_2+k_{-3}+k_4}{S_{tot}(k_1+k_3)}, \frac{(k_{-1}+k_2)(k_{-3}+k_4)}{S_{tot}\big(k_1(k_{-3}+k_4)+k_3(k_{-1}+k_2)\big)},\\
 &\frac{h_{-1}+h_2+h_{-3}+h_4}{S_{tot}(h_1+h_3)}, \frac{(h_{-1}+h_2)(h_{-3}+h_4)}{S_{tot}\big(h_1(h_{-3}+h_4)+h_3(h_{-1}+h_2)\big)}\bigg\},
\end{align*}
then the assumptions {\bf A1}, {\bf A2} and {\bf A4} will hold.

Notice that $\dot{y}$ in \eqref{eqn:ex_full} is linear in $y$ when
$\varepsilon=0$, 
so $g_1$ (defined as in \eqref{eqn:f_xz}) is linear in $z$, and hence the equation
for $z$ can be written as: 
\[
\frac{dz}{d\tau} = B(x_0)z, \ \ x_0 \in U,
\]
where the matrix $B(x_0)$ is Hurwitz for every $x_0 \in U$. Therefore, {\bf A3} also holds.

To check {\bf A6}, let us look at the reduced system ($\varepsilon=0$ in \eqref{eqn:ex_full}):
\begin{align}
\label{eqn:app_limit}
\frac{dx_1}{dt}=&-\frac{k_2x_1}{\frac{K_{m1}}{S_{tot}}+\frac{K_{m1}(1-x_1-x_2)}{K_{m2}}+x_1} +\frac{h_4c\frac{K_{m3}(1-x_1-x_2)}{K_{m4}}}{\frac{K_{m3}}{S_{tot}}+\frac{K_{m3}(1-x_1-x_2)}{K_{m4}}+x_2}:=F_1(x_1,x_2) \\
\frac{dx_2}{dt}=&-\frac{h_2cx_2}{\frac{K_{m3}}{S_{tot}}+\frac{K_{m3}(1-x_1-x_2)}{K_{m4}}+x_2} +\frac{k_4\frac{K_{m1}(1-x_1-x_2)}{K_{m2}}}{\frac{K_{m1}}{S_{tot}}+\frac{K_{m1}(1-x_1-x_2)}{K_{m2}}+x_1}:=F_2(x_1,x_2). \notag 
\end{align}
It is easy to see that $F_1$ is strictly decreasing in $x_2$, and $F_2$ is strictly decreasing in $x_1$ on 
\[
K_0=\{(x_1,x_2) \,|\, x_1 \geq 0, x_2 \geq 0, x_1+x_2 \leq 1\}.
\]
So, \eqref{eqn:app_limit} is strongly monotone on some open set $W$ containing $K_0$ with respect to the cone
\[
\{(x_1,x_2)\, | \, x_1 \leq 0, x_2 \geq 0\}.
\]
Applying Lemma \ref{pro:ESM_EPD}, the flow of \eqref{eqn:app_limit} has eventually positive derivatives on $K_0$ ($\subset W$), and {\bf A6} is valid.

So the system formulated in the form of \eqref{eqn:ex_full} satisfies all
assumptions {\bf A1} to {\bf A7}. Applying Theorem \ref{thm:main}, we have:

\begin{thm}
\label{thm:ex}
There exist a positive $\varepsilon^*<\varepsilon_0$ such that for each $\varepsilon \in (0,\varepsilon^*)$, the forward trajectory of \eqref{eqn:ex_full} starting from almost every point in $D_\varepsilon$ converges to some equilibrium. 
\end{thm}

In fact, since the reduced system is of dimension two, we know that {\em every} trajectory in $D_\varepsilon$, instead of {\em almost every} trajectory in $D_\varepsilon$, converges to some equilibrium (\cite{Hirsch-Smith}). 

It is worth pointing out that the conclusion we obtained from the above
theorem is only valid for small enough $\varepsilon$; that is, the
concentration of the enzyme should be much smaller than the concentration of the
substrate. Unfortunately, this is not always true in biological systems,
especially when feedbacks are present.
However, if the sum of the Michaelis-Menten
constants and the total concentration of the substrate are much larger
than the concentration of enzyme, a different scaling:
\[
x_1=\frac{[S_0]}{A}, \ \ x_2=\frac{[S_2]}{A}, \ \ \varepsilon'=\frac{E_{tot}}{A}, \ \ t=\tau \varepsilon',
\]
where $A=S_{tot}+K_{m1}+K_{m2}+K_{m3}+K_{m4}$ will allow us to obtain the same convergence result.

\subsection{Another example}
The following example demonstrates the importance of the smallness of
$\varepsilon$. Consider an $m+1$ dimensional system:
\begin{align}
\label{eqn:example}
\frac{dx}{dt}=&\gamma(y_1,\dots,y_m)-\beta(x)  \\
\varepsilon\frac{dy_i}{dt}=&-d_iy_i-\alpha_i(x), \ \ d_i>0, \ \ i=1, \dots, m. \notag
\end{align}
under the following assumptions:
\begin{enumerate}
\item There exists an integer $r >1$ such that the derivatives of $\gamma,
  \beta$, and $\alpha_i$ are of class $C^r_b$ for sufficiently large bounded sets.
\item The function $\beta(x)$ is odd, and it approaches infinity as $x$ approaches infinity.
\item The function $\alpha_i(x)$ ($i=1,\dots, m$) is bounded by positive constant $M_i$ for all $x \in \mathbb R$.
\item The number of roots to the equation 
\[
\gamma(\alpha_1(x), \dots, \alpha_m(x))=\beta(x)
\]
is countable.
\end{enumerate}
We are going to show that on any large enough region, and provided that
$\varepsilon$ is sufficiently small, almost every trajectory converges to an equilibrium. To emphasize the need for small $\varepsilon$, we also show that
when $\varepsilon>1$, limit cycles may appear.

Assumption 4 implies {\bf A7}, and because of the form of \eqref{eqn:example}, {\bf A3} and {\bf A4} follow naturally. {\bf A6} also holds, as every one dimensional system is strongly monotone. For {\bf A5}, we take
\[
D_\varepsilon=\{(x,y) \, | \, |x| \leq a, |y_i| \leq b_i, i =1, \dots, m \},
\]
where $b_i$ is an arbitrary positive number greater than $\frac{M_i}{d_i}$ and $a$ can be any positive number such that 
\[
\beta(a)>N_b:=\max_{|y_i| \leq b_i} \gamma(y_1, \dots, y_m).
\]
Picking such $b_i$ and $a$ assures 
\[
x \frac{dx}{dt}<0, \ \ y_i \frac{d y_i}{d t}<0,
\]
i.e. the vector field points transversely inside on the boundary of $D_\varepsilon$. Let $U$ and $V$ be some bounded open sets such that $D_\varepsilon \subset U \times V$, and assumption 1 holds on $U$ and $V$. Then {\bf A1} and {\bf A2} follow naturally.
 By our main theorem, for sufficiently small $\varepsilon$, the forward trajectory of (\ref{eqn:example}) starting from almost every point in $D_\varepsilon$ converges to some equilibrium.

On the other hand, convergence does not hold for large $\varepsilon$. Let
\[
\beta(x)=\frac{x^3}{3}-x, \ \alpha_1(x)=2\tanh x, \ \ m=1,\  \gamma(y_1)=y_1, \ d_1=1.
\]
It is easy to verify that $(0,0)$ is the only equilibrium, and the Jacobian matrix at $(0,0)$ is 
\[
\left( \begin{array}{cc}
         1 & 1 \\
         -2/\varepsilon & -1/\varepsilon
 \end{array} \right).
\]
When $\varepsilon>1$, the trace of the above matrix is
$1-1/\varepsilon>0$, its determinant is $1/\varepsilon>0$,
so the (only) equilibrium in $D$ is repelling. By the Poincar\'e-Bendixson Theorem, there exists a limit cycle in $D$.

\section{Conclusions}
\label{conclusion}

Singular perturbation techniques are routinely used in the analysis of
biological systems.  The \emph{geometric} approach is a powerful tool for
global analysis, since it permits one to study the behavior for finite
$\varepsilon$ on a manifold in which the dynamics is ``close'' to the slow
dynamics.  Moreover, and most relevant to us, a suitable fibration structure
allows the ``tracking'' of trajectories and hence the lifting to the full
system of the exceptional set of non-convergent trajectories, if the slow
system satisfies the conditions of Hirsch's Theorem.  Using the geometric
approach, we were able to provide a global convergence theorem for singularly
perturbed strongly monotone systems, in a form that makes it applicable to the
study of double futile cycles.

\section*{Acknowledgment}

 We have benefited greatly from correspondence with Christopher Jones and
 Kaspar Nipp about geometric singular perturbation theory.  We also wish to
 thank Alexander van Oudenaarden for questions that triggered much of this
 research, and David Angeli, Thomas Gedeon, and Hal Smith for helpful
 discussions.

This work was supported in part by NSF Grant DMS-0614371.

\end{document}